\documentclass{amsart}
\usepackage{amssymb,latexsym}

\newtheorem{definition}{Definition}
\newtheorem{theorem}{Theorem}
\newtheorem{lemma}{Lemma}

\newtheorem{proposition}{Proposition}

\newtheorem{problem}{Problem}

\newcommand{\K}{{\bf k}}

\newcommand{\Z}{{\mathbb Z_{\leq 0}}}

\newcommand{\B}{{\frak B}}
\newcommand{\A}{\mathcal{A}_{\infty}}

\newcommand{\dgart}{{\texttt{dgart}}}
\newcommand{\art}{{\texttt{art}}}
\newcommand{\Alg}{{\texttt{Alg}}}
\newcommand{\DAlg}{{\texttt{DGAlg}}}
\newcommand{\Grp}{{\texttt{Grp}}}
\newcommand{\GRP}{{\underline{\texttt{Grp}}}}
\newcommand{\SET}{{\texttt{SSet}}}
\newcommand{\Mor}{{\texttt{Mor}}}
\newcommand{\Def}{{\texttt{Def}}}
\newcommand{\g}{{\frak{g}}}
\newcommand{\og}{{\overline{g}}}
\newcommand{\od}{{\overline{d}}}

\newcommand{\m}{{\texttt{m}}}

\newcommand{\li}{\mathcal{L}_\infty}

\newcommand{\DEF}{{\underline{\texttt{Def}}}}
\newcommand{\DEL}{{\underline{\texttt{Del}}}}
\newcommand{\HOM}{{\underline{Hom}}}
\newcommand{\Mod}{{\texttt{dgmod}}}
\newcommand{\Deg}{{\texttt{deg}}}
\newcommand{\al}{{\widehat{\alpha}}}
\newcommand{\be}{{\widehat{\beta}}}
\newcommand{\muu}{{\widehat{\mu}}}
\newcommand{\nuu}{{\widehat{\nu}}}
\newcommand{\phii}{{\widetilde{\phi}}}
\newcommand{\de}{{\widehat{\delta}}}
\newcommand{\D}{{\widehat{d}}}
\newcommand{\chai}{{\widehat{\chi}}}
\newcommand{\la}{{\widehat{\lambda}}}
\newcommand{\ro}{{\widehat{\rho}}}
\newcommand{\Sy}{{\overline{S}(s\frak{g}\tensor R)}}
\newcommand{\ga}{{\widetilde{\gamma}}}
\newcommand{\h}{{\frak{h}}}
\newcommand{\tensor}{{\otimes_{\K}}}
\newcommand{\Tensor}{{\otimes_R}}
\newcommand{\FIBER}{{\texttt{h-fib}}}
\newcommand{\DALG}{{\underline{\texttt{DGAlg}}}}
\newcommand{\SC}{{\underline{\texttt{Cat}}}}
\newcommand{\PI}{{\underline{\pi^*}}}
\newcommand{\GAMMA}{{\underline{\mathbb L\gamma}}}
\newcommand{\Morf}{{\texttt{Mor}^{gr,fl}}}
\newcommand{\Morc}{{\texttt{Mor}^{w,c}}}
\newcommand{\F}{{\underline{F}_R}}
\newcommand{\BAMMA}{{\underline{\Omega\B\gamma}}}

\usepackage{thumbpdf}
\usepackage{hyperref}

\begin{document}
\title{Formal deformations of morphisms of associative algebras}
\author{Dennis V. Borisov}
\date{\today}
\maketitle

\begin{abstract}
The $\li$-algebra controlling simultaneous deformations of a morphism of associative algebras and its domain and codomain is described.
Isomorphism of the cohomology of this $\li$-algebra with the classical construction is shown.
\end{abstract}

\tableofcontents

In this work we describe simultaneous deformations of a morphism of dg associative algebras together with the domain and the codomain. By
associative algebras we mean non-unital ones. Such deformations are described not by a Lie but by a proper $\li$-algebra. M.Gerstenhaber and
S.D.Schack have described the Lie structure on the cohomology of this $\li$-algebra (\cite{GS3},\cite{GS2}). Note that if we keep the codomain
constant (as in \cite{vdL}(5.5)) the $\li$-algebra becomes a Lie algebra. One could get a Lie algebra in the general case too by starting with
the colored operad of morphisms of associative algebras (\cite{Mar1}) and then taking the usual tangent Lie algebra in the colored context.

In section \ref{Begining} we outline the theory of formal homotopical deformations, developed in \cite{Hin1}. We work in characteristic 0 and
hence consider the homotopy theory in the category of dg algebras. However, to capture the higher equivalences (i.e. work with
``$\infty$-groupoids'' instead of the usual ones) we have to take into account not only the localization of the category by quasi-isomorphisms,
but also all the mapping spaces between the objects. The starting point is the simplicial localization (\cite{DK3}) of the category of morphisms
of algebras. Then we use simplicial structure on this category (characteristic 0) and express homotopy deformations as fibers of functors
between simplicial categories. The deformation problem then is to find an $\li$-algebra whose simplicial Deligne groupoid is equivalent to those
fibers.

In section \ref{A} we construct this $\li$-algebra. The starting point is a translation of an associative deformation problem into an $\A$ one
(for the case of a single algebra see \cite{PS}). This allows us to express deformation conditions as an equation on Hochschild cochains of the
initial algebras. This equation is the basis for defining the $\li$-algebra, eventually it becomes the Maurer-Cartan equation. At the end of the
section we describe equivalences between deformations. Contrary to the single algebra case, not all derivations generate infinitesimal
automorphisms, hence the Deligne groupoid is built by defining objects to be solutions of the Maurer-Cartan equation in the whole $\li$-algebra
and morphisms being exponentials of a Lie subalgebra. This algebra is a proper $\li$-algebra and not a dg Lie one. The reason is that the
defining equation of an algebra morphism $A\rightarrow B$ requires composition of elements of $Hom(\underset{n>0}{\bigoplus}A^{\tensor^n},B)$
and $Hom(B^{\tensor^m},B)$, where we take $m$ elements from the former and plug them into one from the latter. This is an operation with $m+1$
inputs, hence even if $B$ has only a binary multiplication, we would have ternary operations in the $\li$-algebra.

In section \ref{Solution} we prove that for the morphisms between non-positively graded algebras, the $\li$-algebra we constructed indeed
describes deformations. There are two things to prove there: firstly the Deligne groupoid of the $\li$-algebra should be equivalent to the
groupoid of deformations, i.e. every deformation should be equivalent to the one described by a solution of the Maurer-Cartan equation, secondly
the mapping spaces between deformations should be weakly equivalent to the mapping spaces between solutions, i.e. the simplicial Deligne
groupoid of the $\li$-algebra should be equivalent to the simplicial groupoid of the deformations.

The proof of the first issue uses the fact that $\A$-morphisms are morphisms between codifferential coalgebras and hence there is a natural
embedding of them into the category of associative morphisms, this embedding is homotopically surjective because $\A$-operad is cofibrant. The
second issue is proved by showing that the bar-cobar construction defines weak equivalences between the mapping spaces.

The condition on grading comes from the fact that almost free non-positively graded algebras are cofibrant, which is in general not true for the
$\mathbb Z$-graded ones. Since we describe deformations as fibers of functors between model categories, objects that are not cofibrant are not
in the correct fiber.

Finally we show that the cohomology of our $\li$-algebra coincides with the one in \cite{GS2},\cite{GS3}.

\textbf{Acknowledgements}
    The author is grateful to V.Hinich, B.Keller, S.Merkulov, J.Stasheff and D.Tamarkin for many helpful remarks, and to
    Yu.I.Manin for his attention to the paper.

\textbf{Notation}
    We fix a field $\K$ of characteristic 0. For a $\K$-module $A$ we denote by $S(A)$ the cofree cocommutative non-counital coalgebra,
    cogenerated by $A$ (cofree in the category of connected coalgebras, see e.g. \cite{LM} page 2150). Similarly by $T(A)$ we denote the
    coassociative non-counital coalgebra, that is cofreely cogenerated
    by $A$, i.e. $T(A):=\underset{n>0}{\bigoplus}A^{\tensor^n}$.
    Working with Hochschild cochains
    $C^\bullet(A,A)$ we denote by $\al$ the coderivation of $T(A)$ generated by $\alpha$, and by $\ga$ (for $\gamma\in C^\bullet(A,B)$) the
    morphism of coalgebras $T(A)\rightarrow T(B)$.

    By simplicial categories we mean categories enriched over $\SET$ (category of simplicial sets). Simplicial groupoids are simplicial
    categories, such that the category of connected components is a groupoid. We underline objects and structures that are simplicial:
    $\GRP,\SC$ are categories of simplicial groupoids and categories. $\Grp$ is the category of groupoids. We treat the category of categories as an
    enriched category.

    Differentials raise degree, $\art,\dgart$ are categories of commutative artinian $\K$-algebras and $\Z$-graded dg $\K$-algebras.
    For $R\in\dgart$, $\Alg(R),\DAlg(R)$
    are categories of  associative algebras and dg associative algebras over $R$, $\Mod(R)$ means cochain complexes over $R$. Given a morphism
    $f:R_1\rightarrow R_2$ we denote by $f^*(A):=A\otimes_{R_1}R_2$ the direct image functor $\DAlg(R_1)\rightarrow \DAlg(R_2)$. The unique
    morphism from $R$ to $\K$ is denoted by $\pi_R$.

    Adding $\Mor$ to a
    category name means taking the category of morphisms in that category. Morphisms between morphisms are pairs of morphisms, making
    up a commutative diagram. When we write an element of a module as an exponent we mean its parity.

\section{Homotopical deformations}\label{Begining}
Let $\gamma:A\rightarrow B\in \Mor\Alg(\K)$, we will denote it simply by $\gamma$. Classically (\cite{GS3}) deformations of $\gamma$ are
described by a functor $\Def^\gamma:\art(\K)\rightarrow\Grp$, defined by
    $$\Def^\gamma(R):=(\Morf\Alg(R),\pi^*_R,\gamma),$$
where the image is the comma category of $\pi^*_R$ over $\gamma$, with $\Morf\Alg(R)$ being the underlying groupoid of the category of morphisms
between flat $R$-algebras. If we define equivalences between functors from $\art(\K)$ to $\Grp$ as natural transformations, that are objectwise
equivalences of groupoids, then $\Def^\gamma$ is equivalent to $(\pi^*)^{-1}(Id_\gamma)$, that maps $R$ to the pre-image in $\Morf\Alg(R)$ under
$\pi^*_R$ of the identity automorphism of $\gamma$. This is an instance of the more general weak equivalences of functors to the category of
simplicial groupoids $\GRP$.

To define homotopical deformations of $\gamma$ we consider it as an object in $\Mor\DAlg(\K)$ and extend the parameter category to $\dgart(\K)$.
Weakly equivalent objects should have weakly equivalent deformations (weak equivalences of morphisms are pairs of quasi-isomorphisms). One could
redefine $\Def^\gamma(R)$ as the pre-image of the identity automorphism of $\gamma$ with respect to the functor between the localizations by
quasi-isomorphisms of $\Mor\DAlg(R)$ and $\Mor\DAlg(\K)$, however this will not capture all of the homotopical structure. Instead we will use
the left derived functor of the localization - simplicial localization (\cite{DK3}). So homotopical deformations are described by a functor from
$\dgart(\K)$ to $\GRP$.

Identifying deformations of weakly equivalent objects we are forced to consider simplicial groupoids up to weak equivalences. In order to be
able to do homotopy theory we use a closed model structure on the category of simplicial categories, and then consider the subcategory of
simplicial groupoids. Recall that by simplicial categories we mean categories enriched over simplicial sets. If $\underline{C}$ is a simplicial
category, we denote by $\pi_0(\underline{C})$ the category with the same set of objects and morphisms being the sets of connected components of
the corresponding spaces of maps in $\underline{C}$, and by $\pi_0(\underline{F})$ of a simplicial functor $\underline{F}$, the corresponding
functor between $\pi_0$ of the categories.
\begin{proposition}\label{Bergner}(\cite{Ber})
    Let $\underline{C_1},\underline{C_2}\in\SC$ be simplicial categories. Call a functor $\underline{F}:\underline{C_1}\rightarrow\underline{C_2}$
    \textbf{a weak equivalence} if $\pi_0(\underline{F}):\pi_0(\underline{C_1})\rightarrow\pi_0(\underline{C_2})$ is an equivalence of categories and
    for any $c_1,c_2\in\underline{C}$ the map of $\HOM(c_1,c_2)$ to $\HOM(\underline{F}(c_1),\underline{F}(c_2))$ is a weak equivalence of
    simplicial sets. Call $\underline{F}$ \textbf{a fibration} if for any $c_1,c_2\in\underline{C}$ the map of $\HOM(c_1,c_2)$ to
    $\HOM(\underline{F}(c_1),\underline{F}(c_2))$ is a fibration of simplicial sets, and for any $c_1\in\underline{C_1}$, $x\in\underline{C_2}$ and any
    $f\in\HOM_0(\underline{F}(c_1),x)$, s.t. $\pi_0(f)$ is invertible, there is an $f'\in\HOM_0(c_1,c_2)$, s.t. $\pi_0(f')$ is
    invertible and $\underline{F}(f')=f$. Weak equivalences and fibrations so defined
    are a part of the structure of a closed model category on $\SC$.
\end{proposition}
Note that the weak equivalences in the closed model structure on $\SC$ described in proposition \ref{Bergner} are different from those described
in \cite{Hin1}(5.1.3), where they are required only to induce weak equivalences on the nerves of the categories of connected components. However
for simplicial groupoids these definitions of weak equivalences coincide since a functor between groupoids is an equivalence if and only if it
induces a weak equivalence of the corresponding nerves.

The category $\dgart(\K)$ has a subcategory of weak equivalences consisting of morphisms, that induce isomorphisms on cohomology. Let
$Hom(\dgart(\K),\GRP)$ be the category of functors from $\dgart(\K)$ to $\GRP$ that map weak equivalences to weak equivalences. Simplicial
groupoid describing homotopical deformations of $\gamma$ will be an object in this category.

Let $\underline{C}$ be a simplicial category. Following \cite{DK1}(6.3) we introduce the underlying simplicial groupoid of $\underline{C}$.
Recall that simplicial groupoid means a simplicial category whose $\pi_0$ is a groupoid.
\begin{definition}
    \textbf{The homotopy groupoid} $\underline{C}^{gr}$ of $\underline{C}$ is the maximal simplicial subcategory of $\underline{C}$, which is a
    simplicial groupoid, i.e. objects of $\underline{C}^{gr}$ are those of $\underline{C}$, and for any two of them $\HOM_{\underline{C}^{gr}}
    (c_1,c_2)$ consists of the connected components of $\HOM_{\underline{C}}(c_1,c_2)$, whose classes in $\pi_0(\underline{C})$ are invertible.
\end{definition}

Let $L\Mor\DAlg(R)$ be the simplicial localization (\cite{DK3}) of $\Mor\DAlg(R)$ with respect to quasi-isomorphisms. Any map $f:R_1\rightarrow
R_2$ in $\dgart(\K)$ induces a functor $f^*$ from the category of morphisms of algebras over $R_1$ to those over $R_2$. In general this functor
does not preserve weak equivalences, but it does so for weak equivalences between cofibrant objects (\cite{Hin3}(3.3)). Also we have a
functorial cofibrant replacements for objects in $\Mor\DAlg(R)$ (\cite{Hin3}(7.4.2)), combining $f^*$ with these replacements we get a functor
that preserves all weak equivalences and hence induces a functor $Lf^*$ between the corresponding simplicial localizations. If $f$ is a weak
equivalence itself, $f^*$ is a part of a Quillen equivalence, and hence $Lf^*$ is a weak equivalence between simplicial categories
(\cite{DK1}(3.6)).
\begin{definition}\label{FirstDefinition}
    \textbf{Homotopical deformations} of $\gamma$ are described by the functor $L\Def^\gamma\in Hom(\dgart(\K),\GRP)$, defined by
        $$L\Def^\gamma(R):=\FIBER_{\underline{\gamma}}(L\pi^*_R:(L\Mor\DAlg(R))^{gr}\rightarrow (L\Mor\DAlg(\K))^{gr}),$$
    where $\FIBER_{\underline{\gamma}}$ stands for homotopy fiber at $\underline{\gamma}$, and $\underline{\gamma}$ is the image in
    $L\Mor\DAlg(R)$ of the final object in $\SC$, given by $\gamma$ and the trivial simplicial set consisting of $Id_\gamma$.
\end{definition}

Definition \ref{FirstDefinition} uses simplicial localization, which is hardly computable. In order to describe homotopical deformations
effectively we will use the simplicial model structure on $\Mor\DAlg(R)$. Let $\DALG(R)$ be the simplicial model category (\cite{Hin3}(4.8))
with the same objects as $\DAlg(R)$ and the spaces of maps defined by
    $$\HOM_n(A,B):=Hom(A,B_n),\quad A,B\in\DAlg(R),\quad B_n:=B\Tensor(R\tensor\Omega_n),$$
where $\Omega_n$ is the algebra of polynomial forms on the $n$-simplex (\cite{Bou}(1.)), i.e. it is the free commutative unital $\K$-algebra
generated by $\{t^n_0,...,t^n_n,dt^n_0,...,dt^n_n\}$ with $\Deg(t^n_i)=0$, $\Deg(dt^n_i)=1$, and the differential defined by $d(t^n_i):=dt^n_i$.
The set $\{\Omega_n\}_{n\geq 0}$ is a simplicial algebra, with face and degeneracy maps $\delta_i^n:\Omega_n\rightarrow\Omega_{n-1}$,
$\sigma_i^n:\Omega_n\rightarrow\Omega_{n+1}$ given by
    $$\delta^n_i(t^n_j)=\{t^{n-1}_{j-1}\text{ if }i<j,\text{ 0 if }i=j,\text{ }t^{n-1}_j\text{ else}\},$$
    $$\sigma^n_i(t^n_j)=\{t^{n+1}_{j+1}\text{ if }i<j,\text{ }t^{n+1}_j+t^{n+1}_{j+1}\text{ if }i=j,\text{ }t^{n+1}_j\text{ else}\}.$$

Since category of morphisms of algebras is the category of diagrams in a cofibrantly generated simplicial model category, it also has a
simplicial model structure (\cite{Hir}(11.7) and \cite{ShSh} for the cofibrant generation), defined as follows. Let $\gamma:A\rightarrow B$,
$\gamma':A'\rightarrow B'$ be two objects of $\Mor\DAlg(R)$. A morphism from $\gamma$ to $\gamma'$ is a pair $\phi:A\rightarrow A'$,
$\psi:B\rightarrow B'$. It is a weak equivalence (fibration) if both $\phi$ and $\psi$ are. Let $\Mor\DALG(R)$ be the simplicial model category
with the same objects as $\Mor\DAlg(R)$ and the mapping spaces defined by
    $$\HOM_n(\gamma,\gamma'):=Hom(\gamma,\gamma'_n),\quad\gamma'_n:=(\gamma'\Tensor Id):A'_n\rightarrow B'_n.$$

Classical definition was equivalent to a fiber in the subcategory of flat $R$-algebras and isomorphisms. In homotopy theory flatness condition
is expressed by cofibrant objects and isomorphisms by weak equivalences between cofibrant objects. Let $\Morc\DALG$ be the simplicial
subcategory of $\Mor\DALG(R)$ consisting of cofibrant objects, and the morphisms being
    $$\HOM_n^{w,c}(\gamma,\gamma'):=\{\text{weak equivalences from }\gamma\text{ to }\gamma'_n\}.$$
Let $\PI$ be the simplicial extension of $\pi^*$. Let $\mathbb L\gamma$ be a cofibrant replacement of $\gamma$. The trivial simplicial set
$\{(Id_{\mathbb L\gamma}\Tensor Id):\mathbb L\gamma\rightarrow\mathbb L\gamma_n\}_{n\geq 0}$ together with $\mathbb L\gamma$ itself is an image
in $\Morc\DALG(R)$ of the final object on $\SC$. We will denote it by $\GAMMA$.
\begin{definition}\label{SecondDefinition}(\cite{Hin1})
    \textbf{Homotopical deformations} of $\gamma$ are described by the functor $\DEF^\gamma\in Hom(\dgart(\K),\GRP)$,
    defined by
        $$\DEF^\gamma(R):=\FIBER_\GAMMA(\PI_R:\Morc\DALG(R)\rightarrow\Morc\DALG(\K)),$$
    where $\FIBER_\GAMMA$ stands for homotopy fiber at $\GAMMA$.
\end{definition}

As in the case of definition \ref{FirstDefinition}, $\DEF^\gamma$ is an object of $Hom(\dgart(\K),\GRP)$ because for any morphism
$f:R_1\rightarrow R_2$ in $\dgart(\K)$ the functor $f^*:\Mor\DAlg(R_1)\rightarrow\Mor\DAlg(R_2)$ is a left Quillen functor, i.e. maps
cofibrations (weak equivalences between cofibrants) to the like. Hence it induces a functor
$\underline{f^*}:\Morc\DALG(R_1)\rightarrow\Morc\DALG(R_2)$. Finally $\DEF^\gamma$ maps weak equivalences in $\dgart(\K)$ to weak equivalences
because Quillen equivalences induce weak equivalences of function complexes (\cite{Hir}(17.4.16)).

Since $\mathbb L$ is functorial and there is a natural transformation $\mathbb L\rightarrow Id_{\Mor\DAlg(R)}$, it induces a weak
self-equivalence on $L\Mor\DAlg(R)$, hence $L\Def^\gamma$ and $L\Def^{\mathbb L\gamma}$ are weakly equivalent (homotopical invariance of
$\FIBER$). According to \cite{DK2}(4.8) and \cite{DK1}(2.2) $L\Mor\DAlg(R)$ and $\Mor\DALG(R)$ are naturally weakly equivalent. Therefore the
corresponding homotopical groupoids are weakly equivalent, and since the ``2 out of 3'' axiom implies that $\Morc\DALG(R)$ is the homotopical
groupoid on the cofibrant objects, we conclude that $L\Def^{\mathbb L\gamma}$ and $\DEF^\gamma$ are weakly equivalent. So these two definitions
of homotopical deformations coincide.

\begin{lemma}\label{Fiber}
    The fiber of $\PI_R:\Morc\DALG(R)\rightarrow\Morc\DALG(\K)$ at $\GAMMA$ is weakly equivalent to the homotopy fiber.
\end{lemma}
\textbf{Proof:}
    As it was noted before, for simplicial groupoids the notion of weak equivalences from proposition \ref{Bergner} coincides with that of
    \cite{Hin1}(5.). Also fibrations between simplicial groupoids according to proposition \ref{Bergner} are fibrations according to
    \cite{Hin1}, indeed right lifting property with respect to adding an ingoing or an outgoing arrow (\cite{Hin1}(5.1.4)) is held by fibrations
    between simplicial groupoids according to proposition \ref{Bergner}, because those functors map $\HOM_0$'s surjectively. Therefore according
    to \cite{Hin1}(5.3.2) the nerve functor maps weak equivalences and fibrations
    (as in proposition \ref{Bergner}) between simplicial groupoids to the like between simplicial sets.

    $\PI_R:\Morc\DALG(R)\rightarrow\Morc\DALG(\K)$ is a fibration of simplicial groupoids (\cite{Hin1}(4.2.1)). Indeed, every morphism in $\Morc\DALG(\K)$
    has a pre-image in $\Morc\DALG(R)$, for example $R$-linear extension. If $\gamma_1,\gamma_2\in\Morc\DALG(R)$, then
    $\PI:\HOM(\gamma_1,\gamma_2)\rightarrow\HOM(\PI(\gamma_1),\PI(\gamma_2))$ is a fibration because it coincides with the map
    $\HOM(\gamma_1,\gamma_2)\rightarrow\HOM(\gamma_1,\gamma_2\Tensor\K)$ given by the morphism $\gamma_2\rightarrow\gamma_2\Tensor\K$, and this last
    map is a fibration whereas $\gamma_1$ is a cofibrant object (\cite{Hir}(9.3.1)).

    The nerve functor has a left adjoint (\cite{Hin1}(5.3.1)), hence it preserves limits and therefore maps the nerve of the fiber of a functor
    to the fiber of the nerve of that functor. Since $\SET$ is a proper model category, taking fibers of fibrations is equivalent to taking homotopy
    fibers (\cite{Hir}(13.4.6)). Hence we conclude that the canonical map from the homotopy fiber of $\PI_R:\Morc\DALG(R)\rightarrow
    \Morc\DALG(\K)$ to the fiber induces weak equivalence of the nerves, hence it is a weak equivalence of the groupoids themselves according to
    \cite{DK3}(9.6) (we can consider the two fibers to have the same set of objects because the functor between them induces equivalence on
    the underlying categories of connected components).
$\blacksquare$

Lemma \ref{Fiber} reduces the question of computing deformation groupoid to computing a fiber of a functor. As it is explained in \cite{Hin2},
(the nerve of) such functor can be described by the nerve of a dg Lie algebra. Usually one looks for such dg Lie algebra appearing naturally
from the structure to be deformed. In case of deformation of one algebra it would be the dg Lie algebra of derivations of a cofibrant
replacement. Equivalences between solutions of the Maurer-Cartan equation in this dg Lie algebra are represented by infinitesimal automorphisms
of the cofibrant replacement, and when we consider deformations over graded Artin rings, all elements of the dg Lie algebra represent (graded)
automorphisms. However, in case of deformations of morphisms, not all derivations are infinitesimal automorphisms, and the algebra is not Lie
but an $\li$-algebra. In such more general situations we can have an $\li$-algebra $\g$ and an $\li$-subalgebra $\h$, that represents
infinitesimal automorphisms. Definition of the Deligne groupoid of a pair $(\h,\g)$ is given in definition \ref{Deligne}.
\begin{problem}
    Let $\gamma$ be a morphism of dg associative algebras. \textbf{The deformation problem} defined by $\gamma$ is to find a pair of
    $\li$-algebras $(\h,\g)$, s.t. the corresponding simplicial Deligne groupoid $\DEL^\gamma$ is weakly equivalent to $\DEF^\gamma$.
\end{problem}
In case $\gamma$ is a morphism of non-positively graded algebras, solution to this problem is given in section \ref{Solution}. We will do it by
expressing this deformation problem in the language of $\A$-structures. That will allow us to compare the result with the Lie algebra on
cohomology defined in \cite{GS2} and \cite{GS3}.

\section{Deformation of $\A$-structures}\label{A}

\subsection{Definition of $\A$-structures}

We will define an $\A$-structure on a pair of modules as a morphism between two coassociative codifferential coalgebras. This requires a choice
of identification of the tensor algebra on a module with the one on its suspension. We make the following choice (\cite{Pen}(5.)). Let
$M\in\Mod(R)$, then $(sM)^k:=M^{k+1}$. Let $T(M):=\underset{n>0}{\bigoplus}M^{\Tensor^n}$, then we define
        $$s:T(M)\rightarrow T(sM),\quad s(x_1\Tensor...\Tensor x_n):=(-1)^\epsilon sx_1\Tensor...\Tensor sx_n,$$
        $$\epsilon=(n-1)x_1+(n-2)x_2+...+x_{n-1}.$$

Recall that by $\al$ we mean the coderivation of a cofree coalgebra, cogenerated by a map $\alpha$, and by $\ga$ we mean the coalgebra morphism,
cogenerated by $\gamma$.
\begin{definition}\label{DefinitionAstructure}(e.g. \cite{Pen}(5.))
    Let $R\in\dgart(\K)$. Let $A,B\in\Mod(R)$. \textbf{An $\A$-structure} on the pair $(A,B)$ is a triple $(\al,\be,\ga)$, where
    $\al,\be$ are coderivations of degree 1 on the cofree coassociative coalgebras $T(sA),T(sB)$ respectively and $\ga:T(sA)\rightarrow T(sB)$
    is a degree 0 morphism of coalgebras, such that
        $$(\al+\de_A)^2=0,\quad(\be+\de_B)^2=0,\quad\ga\circ(\al+\de_A)-(\be+\de_B)\circ\ga=0,$$
    where $\de_A,\de_B$ are the coderivations cogenerated by differentials on $A,B$.
\end{definition}
Since $T(sA),T(sB)$ are cofree coalgebras, $\al,\be,\ga$ are determined by their corestrictions to cogenerators
    $$\{\alpha_n:(sA)^{\Tensor^n}\rightarrow sA\}_{n>0},\quad\{\beta_n:(sB)^{\Tensor^n}\rightarrow sB\}_{n>0},\quad\{\gamma_n:(sA)^{\Tensor^n}\rightarrow
    sB\}_{n>0}.$$
Using identification $s$ we can translate $\alpha,\beta,\gamma$ to tensor operations on $A,B$, namely
    $$\muu:=s^{-1}\circ\al\circ s,\quad\nuu:=s^{-1}\circ\be\circ s,\quad\phii:=s^{-1}\circ\ga\circ s.$$
Applying this to the corestrictions we get the sequences of Hochschild cochains $\{\mu_n\in C^n(A,A)\}_{n>0}$, $\{\nu_n\in C^n(B,B)\}_{n>0}$,
$\{\phi_n\in C^n(A,B)\}_{n>0}$, such that in the case $A,B$ have trivial differentials
    $$\Deg(\mu_n)=\Deg(\nu_n)=2-n,\quad\Deg(\phi_n)=1-n,$$
    $$\underset{\underset{0\leq i\leq n-k}{k+l=n+1}}{\Sigma}(-1)^{\epsilon_1}\mu_l(a_1...\mu_k(a_{i+1}...a_{i+k})...a_n)=0,\quad
    \epsilon_1=i(k-1)+k(n-k+a_1+...+a_i),$$
    $$\underset{\underset{0\leq i\leq n-k}{k+l=n+1}}{\Sigma}(-1)^{\epsilon_2}\nu_l(b_1...\nu_k(b_{i+1}...b_{i+k})...b_n)=0,\quad
    \epsilon_2=i(k-1)+k(n-k+b_1+...+b_i),$$
    $$\underset{1\leq i\leq n-m}{\underset{m+k=n+1}{\Sigma}}(-1)^{\epsilon_3}\phi_k(a_1...\mu_m(a_{i+1}...a_{i+m})...a_n)=
    \qquad\qquad\qquad\qquad\qquad\qquad\qquad\qquad$$
    $$\qquad\qquad\qquad\qquad\qquad\qquad
    =\underset{i_1+...+i_r=n}{\underset{1\leq r\leq n}{\Sigma}}(-1)^{\epsilon_4}\nu_r(\phi_{i_1}(a_1...a_{i_1})...\phi_{i_r}(a_{n-i_r}...a_n)),$$
    $$\epsilon_3=i(m-1)+m(n-m)+m(a_1+...+a_i),\quad\epsilon_4=\underset{1\leq t<r}{\Sigma}\phi_{i_{t+1}}\underset{s=1}
    {\overset{i_1+...+i_t}{\Sigma}}a_s+\underset{1\leq t<r}{\Sigma}(r-t)\phi_t,$$
this is the usual definition of $\A$-algebras and $\A$-morphisms (e.g. \cite{Kel}(3.1,3.4)).

\subsection{$\li$-algebra on the cochain complex}

The right hand side of the last equation contains compositions of the type $\nu_n(\phi_{i_1}\Tensor...\Tensor\phi_{i_n})$. It is these
compositions that make the Hochschild complex of an $\A$-structure not a Lie algebra but an $\li$-algebra.
\begin{definition}(e.g. \cite{Pen}(6.))
    Let $R\in\dgart$, $M\in\Mod(R)$. The structure of an \textbf{$\li$-algebra} on $M$ is given by a degree 1 coderivation $\D$ on the cofree
    cocommutative non-counital coalgebra (cofree in the category of connected coalgebras, see e.g. \cite{LM} page 2150) cogenerated by
    the suspension $sM$
    of $M$, such that $(\D+\de_M)^2=0$, where $\de_M$ is the coderivation cogenerated by the differential on $M$.
\end{definition}
Since a cocommutative coalgebra is also coassociative there is a canonical embedding of it into $T(sM)$, given by the universal property of
$T(sM)$. Let $S(sM)$ denote this coassociative sub-coalgebra of $T(sM)$. Explicitly
    $$S(sM)\cap(sM)^{\Tensor^n}=<\underset{\sigma}{\Sigma}(-1)^{\epsilon(\sigma;sx_1,...,sx_n)}sx_{\sigma(1)}\Tensor...\Tensor sx_{\sigma(n)}>,$$
where the r.h.s. is the $R$-submodule generated by symmetrizations of all elements of $(sM)^{\Tensor^n}$, and $\epsilon(\sigma;sx_1,...,sx_n)$
is the usual sign of a permutation. The following lemma, generalizing the pre-Lie algebra technique (\cite{GS1}(10.1)), is obvious
(\cite{vdL}(3.8)).
\begin{lemma}\label{RedefinitionLemma}
    Let $M\in\Mod(R)$ and let $\D$ be a degree 1 coderivation on $T(sM)$, that commutes with $\de_M$.
    If $\D^2$ vanishes on $S(sM)$, it defines the structure of an $\li$-algebra on $M$.
\end{lemma}
We will denote $\li$-algebra $(S(sM),\D)$ by $(M,\D)$. From now on we fix $A,B\in\Mod(\K)$. The $\li$-algebra, that describes deformations of
$\A$-structures on the pair $(A,B)$ is built on
    $$\g:=Hom(T(sA),sA)\oplus Hom(T(sB),sB)\oplus s^{-1}Hom(T(sA),sB).$$
Differentials on $A,B$ induce codifferentials on $T(sA),T(sB)$, and we consider $\g$ together with the induced differential. We denote an
element $g$ of $\g$ by $\alpha+\beta+s^{-1}\gamma$. Then $\al$, $\be$ denote the coderivations on $T(sA)$, $T(sB)$, cogenerated by $\alpha$,
$\beta$, and $\ga$ denotes the morphism of coalgebras, cogenerated by $\gamma$.
\begin{definition}\label{DefinitionLalgebra}
    Let the coderivation $\D$ on $T(s\g)$ be defined by its corestriction to cogenerators as follows
        $$d:=\chi_A+\chi_B+\lambda+\underset{m\geq 1}{\Sigma}\rho_m.$$
        $$\chi_A(s\alpha_1\tensor s\alpha_2):=(-1)^{\alpha_1} s(\alpha_1\circ\al_2),\quad
          \chi_B(s\beta_1\tensor s\beta_2):=(-1)^{\beta_1}s(\beta_1\circ\be_2),$$
        $$\lambda(\gamma\tensor s\alpha):=-(-1)^{\gamma}\gamma\circ\al,\quad
          \rho_m(s\beta\tensor\gamma_1\tensor,...,\tensor\gamma_m):=\beta(\gamma_1\tensor...\tensor\gamma_m),$$
    where $m\geq 1$ and $(\gamma_1\tensor...\tensor\gamma_m)$ is the
    linear map $T(sA)\rightarrow (sB)^{\tensor^m}$ with the sign, given by the Koszul sign rule. On the rest of $T(s\g)$ the maps
    $\chi_A,\chi_B,\lambda,\{\rho_m\}_{m\geq 1}$ are defined to be zero.
\end{definition}
\begin{proposition}\label{AlgebraOnCochains}
    Degree of $\D$ is 1 and its square vanishes on $S(s\g)$.
\end{proposition}
\textbf{Proof:}
    Composition of operations has degree zero and the l.h.s. of the defining equations of $\chi$, $\lambda$, $\rho$ have 1 more application of the
    suspension than the corresponding r.h.s., therefore $d$ and hence $\D$ have degree 1.

    Since $\Deg(\D)=1$, $\D^2$ is a coderivation. Therefore it is enough to check that its corestriction to the cogenerators of $S(s\g)$ vanishes.
    From the definition it follows that this corestriction is
        $$\chi_A\circ\chai_A+\chi_B\circ\chai_B+\lambda\circ\chai_A+\lambda\circ\la+\lambda\circ\ro_k+\rho_k\circ\la+\rho_k\circ\chai_B+
        \rho_l\circ\ro_m,$$
    the rest of the compositions vanish identically on $T(s\g)$. We will prove that the above sum is zero on $S(s\g)$ by dividing it into 4
    summands.
\begin{itemize}
\item[1.]
    $\underline{\chi_A\circ\chai_A+\chi_B\circ\chai_B}$
        $$(\chi_A\circ\chai_A)(s\alpha_1\tensor s\alpha_2\tensor s\alpha_3)=(-1)^{\alpha_1+(\alpha_1+\alpha_2)}s((\alpha_1\circ\al_2)\circ\al_3)+$$
        $$+(-1)^{(\alpha_1+1)+\alpha_2+\alpha_1}s(\alpha_1\circ(\al_2\circ\al_3)).$$
    Looking at the signs one sees that the part that is not
    zero is when $\alpha_3$ is not "plugged" into $\alpha_2$.
    But that is taken care of by interchanging $\alpha_2$ and
    $\alpha_3$ (with the correct signs). So
    $\chi_A\circ\chai_A$ vanishes on the symmetrization of
    $s\alpha_1\tensor s\alpha_2\tensor s\alpha_3$ and therefore on
    the whole of $S(s\g)$. The same for $\chi_B\circ\chai_B$.
\item[2.]
    $\underline{\lambda\circ\chai_A+\lambda\circ\la}$
        $$(\lambda\circ\chai_A)(\gamma\tensor s\alpha_1\tensor s\alpha_2)=(-1)^{\gamma+\alpha_1+\gamma+1}
        \gamma\circ(\al_1\circ\al_2),$$
        $$(\lambda\circ\la)(\gamma\tensor s\alpha_1\tensor s\alpha_2)=(-1)^{\gamma+1+(\gamma+\alpha_1)+1}
        (\gamma\circ\al_1)\circ\al_2.$$
    As in the step 1. the part that is not zero vanishes after
    symmetrization, i.e. $\lambda\circ\chai_A+\lambda\circ\la$
    is zero on
    $\gamma\tensor s\alpha_1\tensor s\alpha_2+(-1)^{(\alpha_1+1)(\alpha_2+1)}\gamma\tensor s\alpha_2\otimes s\alpha_1$.
    On the rest of the permutations of $\gamma\tensor s\alpha_1\tensor s\alpha_2$ it vanishes by definition of the maps
    $\chi_A,\lambda$.
\item[3.]
    $\underline{\lambda\circ\ro_m+\rho_m\circ\la}$
        $$(\lambda\circ\ro_m)(s\beta\tensor\gamma_1\tensor...\tensor\gamma_m\otimes s\alpha)=(-1)^{r}
        (\beta\circ(\gamma_1\tensor...\tensor\gamma_m))\circ\al,$$
        $$(\lambda\circ\ro_m)(s\beta\tensor\gamma_1\tensor...\tensor\gamma_i\tensor s\alpha\tensor\gamma_{i+1}\tensor...
        \tensor\gamma_m)=0\quad if\quad i<m,$$
        $$(\rho_m\circ\la)((-1)^{(\alpha+1)(\gamma_{i+1}+...+\gamma_m)}
        s\beta\tensor\gamma_1\tensor...\tensor\gamma_i\tensor s\alpha\tensor\gamma_{i+1}\tensor...
        \tensor\gamma_m)=$$
        $$=(-1)^{t}\beta(\gamma_1\tensor...\tensor(\gamma_i\circ\al)\tensor...\tensor\gamma_m),$$
    where
        $$r=\beta+\gamma_1+...+\gamma_m+1,$$
        $$t=\beta+1+\gamma_1+...+\gamma_{i-1}+\gamma_i+1+(\alpha+1)(\gamma_{i+1}+...+\gamma_m).$$
    The difference is $\alpha(\gamma_{i+1}+...+\gamma_m)+1$, but
    $(-1)^{\alpha(\gamma_{i+1}+...+\gamma_m)}$ is exactly the sign
    that appears when we apply first $\al$ and then
    $\gamma_1\tensor...\tensor\gamma_m$ or
    $\gamma_1\tensor...\tensor(\gamma_i\circ\al)\tensor...\tensor\gamma_m$
    directly. It means that
    $\lambda\circ\ro_m+\rho_m\circ\la$ vanishes on the
    symmetrization of $s\beta\tensor\gamma_1\tensor...\tensor\gamma_m\otimes s\alpha$.
\item[4.]
    $\underline{\rho_m\circ\chai_B+\rho_l\circ\ro_n}$
        $$(\rho_m\circ\chai_B)(s\beta_1\tensor s\beta_2\tensor\gamma_1\tensor...\tensor\gamma_m)=(-1)^{\beta_1}(\beta_1\circ\be_2)
        (\gamma_1\tensor...\tensor\gamma_m),$$
        $$(\rho_m\circ\chai_B)(s\beta_1\tensor\gamma_1\tensor...\tensor\gamma_i\tensor s\beta_2\tensor...\tensor\gamma_m)=0\quad
        if\quad i>0,$$
        $$(\rho_l\circ\ro_n)((-1)^{(\beta_2+1)(\gamma_1+...+\gamma_i)}
        s\beta_1\tensor\gamma_1\tensor...\tensor\gamma_i\tensor s\beta_2\tensor...\tensor\gamma_m)=$$
        $$=(-1)^{r}\beta_1(\gamma_1\tensor...\gamma_i\tensor(\beta_2(\gamma_{i+1}\tensor...\tensor\gamma_{i+n}))\tensor...\tensor\gamma_m),$$
    where
        $$r=(\beta_2+1)(\gamma_1+...+\gamma_i)+\beta_1+1+\gamma_1+...+\gamma_i,\quad l+n=m+1.$$
    By definition of $\be_2$, application of $\rho_n\circ\chai_B$ as in the first
    formula is a sum and its summands are exactly the results of
    application of $\rho_l\circ\ro_m$ as in the last formula.
    Then if the signs are opposite we would conclude that
    $\rho_m\circ\chai_B+\rho_l\circ\ro_n$ vanishes on the
    symmetrization of
    $s\beta_1\tensor s\beta_2\tensor\gamma_1\tensor...\tensor\gamma_m$
    (and therefore on the whole of $S(s\g)$). As it is written,
    the difference in signs is
        $$(-1)^{\beta_2(\gamma_1+...+\gamma_i)+1},$$
    but this is exactly opposite to the difference in signs that appears if we
    apply first $\gamma_1\tensor...\tensor\gamma_m$ and then
    $\beta_2$ at position $i+1$ or we apply
    $$\gamma_1\tensor...\tensor(\beta_2(\gamma_{i+1}\tensor...\tensor\gamma_{i+n}))\tensor...\tensor\gamma_m$$ directly.
\end{itemize}
$\blacksquare$

\subsection{Deligne groupoid}

Proposition \ref{AlgebraOnCochains} together with lemma \ref{RedefinitionLemma} imply that $(\g,d)$ is an $\li$-algebra ($d\in Hom(S(s\g),sg)$).
Indeed operations in the definition of $d$ are defined using composition of maps. The operation of composition is a cocycle for the differential
on $\g$, hence $d$ commutes with that differential. We will denote the part of $d$ lying in $Hom((s\g)^{\tensor^m},s\g)$ by $d_m$.

Let $R\in\dgart(\K)$, let $(\g\tensor R,d)$ be the $R$-linear extension of $(\g,d)$, it is an $\li$-algebra over $R$. We will denote the
corestriction of $d$ to $s\g\tensor R$ again by $\underset{m>0}{\Sigma}d_m$.

In order to describe solutions of the Maurer-Cartan equation as points in a coalgebra we have to complete $S(s\g\tensor R)$. Define
$\Sy:=\underset{m>0}{\Sigma}(S(s\g(R))\cap (s\g\tensor R)^{\Tensor^m})$, it is a cocommutative coalgebra. A point in $\Sy$ is an element
$\og_R$, s.t. $\Delta(\og_R)=\og_R\Tensor\og_R$, where $\Delta$ is the comultiplication, this name comes from considering $\li$-algebras as
formal dg manifolds (\cite{Kon}).

Clearly $C^\bullet(A\tensor R,A\tensor R)\oplus C^\bullet(B\tensor R,B\tensor R)\oplus C^\bullet(A\tensor R,B\tensor R)$ is graded by the number
of the arguments of a cochain, $\g\tensor R$ is the completion with respect to the associated filtration, and $d$ is continuous with respect to
this filtration (because it is given by composition of the cochains). Therefore if we denote $\od:=\underset{m>0}{\Sigma}d_m$, it is well
defined as a function on $\Sy$ with values in $s\g\tensor R$, indeed values of $d_m$, for $m>2$, have at least $m$ arguments.
\begin{definition}
    \textbf{A solution of the Maurer-Cartan equation (MCE)} in $(\g\tensor R,d)$ is a degree 0 element $sg_R$, such that the following equation
    holds
        $$\delta_R(sg_R)+\underset{m>0}{\Sigma}d_m((sg_R)^{\Tensor^m})=0,$$
    where $\delta_R$ is the differential on $\g\tensor R$.
\end{definition}
Note that in \cite{Kon}(4.3) in the definition of the Maurer-Cartan equation for $\li$-algebras there is a coefficient $\frac{1}{m!}$ before
$d_m$. We do not have them here because the canonical embedding of the cofree cocommutative coalgebra cogenerated by $s\g\tensor R$ into the
cofree coassociative one maps $(sg_R)^{\odot_R^m}$ to $m!(sg_R)^{\Tensor^m}$, when $g_R$ is odd.

Since $\Sy$ is the completion of a cofree coalgebra, every degree 0 point is $\og_R=\underset{m>0}{\Sigma}(sg_R)^{\Tensor^m}$, for some $sg_R\in
s\g\tensor R$ of degree 0, and for $sg_R$ to be a solution of MCE is equivalent to $\og_R$ being a cocycle for $\od+\delta_R$. We will sometimes
represent solutions of MCE by the corresponding points.

\begin{proposition}\label{BijectionNoReduction}
    There is a bijection between the set of $\A$-structures on the pair $A\tensor R,B\tensor R$ and the set of solutions of MCE in $(\g\tensor R,d)$.
\end{proposition}
\textbf{Proof:}
    Let $\og_R=\underset{m>0}{\Sigma}(sg_R)^{\Tensor^m}$ be a solution of MCE. Then $g_R=\alpha_R+\beta_R+s^{-1}\gamma_R$,
    where $\alpha_R\in Hom(T(sA\tensor R),sA\tensor R)$, $\beta_R\in Hom(T(sB\tensor R),sB\tensor R)$ of degree 1,
    $\gamma_R\in Hom(T(sA\tensor R),sB\tensor R)$ of degree 0, and they
    satisfy $(\al_R+\de_A)^2=0$, $(\be_R+\de_B)^2=0$, $\gamma_R\circ(\al_R+\de_A)-(\beta_R+\de_B)\circ\ga_R=0$ (recall that $\al,\be$ are
    the coderivations, $\ga$ is the coalgebra  morphism, cogenerated by $\alpha,\beta,\gamma$ respectively), where $\de_A$ is the coderivation,
    cogenerated by the differential on $sA\tensor R$. So they comprise an $\A$-structure
    on $A\tensor R,B\tensor R$. Conversely, any such $\A$-structure consists of an element of
        $$Hom(T(sA\tensor R),sA\tensor R)\oplus Hom(T(sB\tensor R),sB\tensor R)\oplus s^{-1}Hom(T(sA\tensor R),sB\tensor R),$$
    and the equations of definition \ref{DefinitionAstructure} translate to the MCE.
$\blacksquare$

Solutions of MCE in $(\g\tensor R,d)$ represent all of $\A$-structures on $A\tensor R,B\tensor R$. We are interested in those, whose reduction
modulo $\m_R$ is the given one $\gamma$ ($\m_R$ stands for the maximal ideal in $R$). Such $\gamma$ is represented by a solution $g$ of MCE in
$(\g,d)$, and a general procedure associates to it a new $\li$-algebra $(\g,d^g)$, that we will use to represent $\A$-structures with the
correct reduction modulo $\m_R$. This is the $\li$-version of the usual technique in dg Lie algebra: changing the differential by adding to it a
bracket with an odd cocycle. The following lemma describes how a solution of MCE can be used to deform the $\li$-algebra.
\begin{lemma}\label{DeformationLstructure}
    Let $(M,d=\underset{n>0}{\Sigma}d_n)$ be an $\li$-algebra ($M\in\Mod(R)$ for some $R\in\dgart(\K)$). Let $x\in M$ be a solution of
    MCE in $(M,d)$. Define
        $$d^x:T(sM)\rightarrow sM,\quad d^x_n(sy_1\Tensor...\Tensor sy_n):=\underset{m\geq 0}{\Sigma}d_{m+n}Sh((sx)^{\Tensor^m},
        sy_1\Tensor...\Tensor sy_n),$$
    where $Sh$ denotes all shuffles of $(sx)^{\Tensor^m}$ in $sy_1\Tensor...\Tensor sy_n$ (there is no sign change because $sx$ is even).
    Then $(M,d^x=\underset{n>0}{\Sigma}d^x_n)$ is an $\li$-algebra.
\end{lemma}
\textbf{Proof:}
    Since tensoring with $(sx)^{\Tensor^m}$ has degree 0 and $d$ is of degree 1, the composition has degree 1.

    Consider application of $(d^x)^2$ to $\underset{\sigma\in S_n}{\Sigma}(-1)^{\epsilon(\sigma;sy_1,...,sy_n)}sy_{\sigma(1)}\Tensor...\Tensor
    sy_{\sigma(n)}$. The result is the sum
        $$\underset{l+m=n+1}{\Sigma}d^x_l\circ\widehat{d^x_m}(\underset{\sigma\in S_n}{\Sigma}
        (-1)^{\epsilon(\sigma;sy_1,...,sy_n)}sy_{\sigma(1)}\Tensor...\Tensor sy_{\sigma(n)})=$$
        $$=\underset{l+m=n+1}{\Sigma}\underset{\underset{i+j=k}{k\geq 0}}{\Sigma}
        d_{l+i}\circ\widehat{d_{m+j}}(\underset{\sigma\in S_n}{\Sigma}
        (-1)^{\epsilon(\sigma;sy_1,...,sy_n)}Sh((sx)^{\Tensor^k}, sy_{\sigma(1)}\Tensor...\Tensor sy_{\sigma(n)})).$$
    This equality is true because the point in $\overline{S}(sM)$, generated by $sx$, is a cocycle for $\od$.
    For each $k\geq 0$ we have
        $$\underset{\underset{i+j=k}{l+m=n+1}}{\Sigma}d_{l+i}\circ\widehat{d_{m+j}}(\underset{\sigma\in S_n}{\Sigma}
        (-1)^{\epsilon(\sigma;sy_1,...,sy_n)}Sh((sx)^{\Tensor^k}, sy_{\sigma(1)}\Tensor...\Tensor sy_{\sigma(n)}))=$$
        $$=\underset{p+q=n+k+1}{\Sigma}\frac{1}{k!}\underset{\sigma\in S_{n+k}}{\Sigma}(-1)^{\epsilon(\sigma;sx,...,sx,sy_1,...,sy_n)}
        d_p\circ\widehat{d_q}(\sigma((sx)^{\Tensor^k}\Tensor sy_1\Tensor...\Tensor sy_n)),$$
    where $\frac{1}{k!}$ appears because $sx$ is even and interchanging $sx$'s does not affect the sign of the
    permutation, whereas on the l.h.s. of the equation the permutations of $sx$'s are absent. The r.h.s. is
    obviously $\frac{1}{k!}$ times $\widehat{d}^2$, applied to the symmetrization of $(sx)^{\Tensor^k}\Tensor
    sy_1\Tensor...\Tensor sy_n$, therefore for each $k\geq 0$ the corresponding sum is 0.
$\blacksquare$

Now we fix an $\A$-structure $A\overset{\gamma}{\rightarrow}B$ (we will denote it simply by $\gamma$). Let $g$ be the corresponding solution of
MCE in $(\g,d)$ (proposition \ref{BijectionNoReduction}), we extend $R$-linearly the $\li$-structure from $(\g,d^g)$ to $(\g\tensor R,d^g)$ and
then consider the $\li$-subalgebra $(\g\tensor\m_R,d^g)$, it is an $\li$-algebra in $\Mod(R)$.
\begin{proposition}
    Let $g$ be the solution of MCE in $(\g,d)$, that corresponds to $\gamma$.
    There is a bijection between the set of solutions of MCE in $(\g\tensor\m_R,d^g)$ and the set of $\A$-structures on $(A\tensor R,B\tensor R)$,
    whose reduction modulo $\m_R$ is $\gamma$.
\end{proposition}
\textbf{Proof:}
    Let $g'$ be a solution of MCE in $(\g\tensor\m_R,d^g)$, then $g+g'\in\g\tensor R$ and we have
        $$d_n((s(g+g'))^{\Tensor^n})=d_n((sg)^{\tensor^n})+\underset{0\leq
        l<n}{\Sigma}d_n(Sh(((sg)^{\tensor^l}\tensor R),(sg')^{\Tensor^{n-l}})).$$
    i.e.
        $$\od(\underset{n>0}{\Sigma}(sg+sg')^{\Tensor^n})=\od(\underset{n>0}{\Sigma}(sg)^{\tensor^n})+
        \overline{d^g}(\underset{n>0}{\Sigma}(sg')^{\Tensor^n}).$$
    Therefore $g+g'$ is a solution of MCE in $(\g\tensor R,d)$. Let $\gamma+\gamma'$ be the $\A$-structure that corresponds to $g+g'$
    (proposition \ref{BijectionNoReduction}). The cochains that generate this
    $\A$-structure are given by $g+g'$, and since $g'\in\g\tensor\m_R$, reduction modulo $\m_R$ of $\gamma+\gamma'$ is obviously generated by $g$.

    Conversely, let $\gamma+\gamma'$ be an $\A$-structure whose reduction modulo $\m_R$ is $\gamma$. Again by proposition
    \ref{BijectionNoReduction} there is a corresponding solution $g+g'$ of MCE in $(\g\tensor R,d)$. Reduction modulo $\m_R$ of $g+g'$
    has to be $g$ and hence $g'\in\g\tensor\m_R$. Using identification of MCE with the defining equations of $\A$-structures
    (proposition \ref{BijectionNoReduction}) we conclude that $g'$ is a solution of MCE in $(\g\tensor\m_R,d^g)$.
$\blacksquare$

In case of a deformation of one algebra (say $A$), $C^\bullet(A,A)$ is a dg Lie algebra, and given a dg Artin algebra $R$, solutions of MCE in
$C^\bullet(A,A)\tensor \m_R$ are equivalent iff the corresponding structures of an $\A$-algebra on $A\tensor R$ are connected by an invertible
$\A$-morphism, whose reduction modulo $\m_R$ is the identity automorphism. Hence on the set of solutions acts the group $(C^\bullet(A,A)\tensor
\m_R)_0$ (with the Campbell-Hausdorff multiplication).

In case of a deformation of the morphism $A\rightarrow B$, equivalences between $\A$-structures are given by pairs of $\A$-morphisms
$A\rightarrow A$ and $B\rightarrow B$, therefore elements of the subspace $\g\cap s^{-1}Hom(T(sA),sB)$ do not represent infinitesimal
automorphisms, instead these are given by the following subspace:
    $$\h:=\g\cap (Hom(T(sA),sA)\oplus Hom(T(sB),sB)).$$
From definition \ref{DefinitionLalgebra} it follows that on $\h$ all ternary and higher operations vanish (indeed the operations that involve 3
and more elements require elements of $s^{-1}Hom(T(sA),sB)$ as inputs), hence, if we forget the differential, $\h\tensor\m_R$ is a nilpotent Lie
algebra. Let $H_R$ be the group, defined on the degree 0 part of $\h\tensor\m_R$, with the Campbell-Hausdorff multiplication. We are going to
define an action of $H_R$ on the set of solutions of MCE in $(\g\tensor\m_R,d^g)$, where $g$ is the solution of MCE in $(\g,d)$, corresponding
to $\gamma$. This action is defined through the adjoint representation of $\g$ on itself (as an $\li$-algebra). This representation comes from
the structure of a left $\li$-module of $\g$ over itself.
\begin{lemma}\label{Equivalences}
    Let $g'$ be a solution of MCE in $(\g\tensor\m_R,d^g)$. For an $h\in H_R$ define
        $$ad_hg':=d^{g+g'}_1(sh)+\delta_R(sh)=\underset{n\geq 0}{\Sigma}d_{n+1}(Sh((sg+sg')^{\Tensor^n},sh))+\delta_R(sh),$$
    where $\delta_R$ is the differential on $\g\tensor\m_R$, then $\underset{k\geq
    0}{\Sigma}\frac{1}{k!}(ad_h)^kg'$ is also a solution of MCE in $(\g\tensor\m_R,d^g)$ and this defines an action of $H_R$ on the set of
    solutions.
\end{lemma}
\textbf{Proof:}
    Consider the groups $G_1,G_2$ that are the subgroups of $Aut(T(sA\tensor R))$, $Aut(T(sB\tensor R))$, consisting of the elements
    whose reduction modulo $\m_R$ are identities on $T(sA),T(sB)$. There is an action of $G_1\times G_2$ on the set of solutions of MCE in
    $(\g\tensor\m_R,d^g)$, given by
        $$\al+\al'\mapsto\phi\circ(\al+\al'+\de_R)\circ\phi^{-1},\quad\be+\be'\mapsto\psi\circ(\be+\be'+\de_R)\circ\psi^{-1},
        \quad\widetilde{\gamma+\gamma'}\mapsto\psi\circ(\widetilde{\gamma+\gamma'})\circ\phi^{-1},$$
    where $(\phi,\psi)\in G_1\times G_2$, $g=\alpha+\beta+s^{-1}\gamma$, $g'=\alpha'+\beta'+s^{-1}\gamma'$.
    For any such $\phi,\psi$ there are $\mu,\nu\in\g\tensor\m_R$, s.t. $\phi=e^{\widehat{\mu}},\psi=e^{\widehat{\nu}}$, where
    $\widehat{\mu},\widehat{\nu}$ are the coderivations cogenerated by $\mu,\nu$. We have
        $$\psi\circ(\widetilde{\gamma+\gamma'})\circ\phi^{-1}=e^{\widehat{\nu}}\circ(\widetilde{\gamma+\gamma'})\circ e^{-\widehat{\mu}}=
        \underset{k\geq 0}{\Sigma}\frac{1}{k!}(ad_{\widehat{\mu}+\widehat{\nu}})^k(\widetilde{\gamma+\gamma'}),$$
    where $ad_{\widehat{\mu}+\widehat{\nu}}(\widetilde{\gamma+\gamma'}):=\widehat{\nu}\circ(\widetilde{\gamma+\gamma'})-
    (\widetilde{\gamma+\gamma'})\circ\widehat{\mu}$. Summands on the r.h.s. of the last equation are $\widetilde{\gamma+\gamma'}$-coderivations
    from $T(sA\tensor R)$ to $T(sB\tensor R)$, therefore their sum is determined by its corestriction to the cogenerators of $T(sB\tensor R)$, and this
    is exactly the projection of $ad_{\mu+\nu}g'$ onto $s^{-1}Hom(T(sA),sB)\tensor\m_R$, similarly for $\al,\be$. So $\underset{k\geq
    0}{\Sigma}\frac{1}{k!}(ad_h)^k$ represents the action of $G_1\times G_2$ on the set of solutions and hence it defines an action of $H_R$,
    since Campbell-Hausdorff multiplication represents composition in the corresponding group.
$\blacksquare$

Using lemma \ref{Equivalences} we can represent morphisms between $\A$-structures by elements of the Lie algebra $\h$. However, we are
interested in the whole spaces of maps. Let $\Omega_n$ be the algebra of polynomial forms on the $n$-simplex described in section
\ref{Begining}. Let $g$ be the solution of MCE in $(\g,d)$, corresponding to $\gamma$.
\begin{definition}\label{Deligne}(\cite{Hin1}(3.1))
    \textbf{The simplicial Deligne groupoid} $\DEL^\gamma(R)$ is given by
        $$Obj(\DEL^\gamma(R)):=\{\text{ the set of solutions of MCE in }(\g\tensor\m_R,d^g)\},$$
        $$\HOM_n(g_1,g_2):=\{h\in\h\tensor\m_R\tensor\Omega_n\text{ s.t. }\underset{k\geq 0}{\Sigma}\frac{1}{k!}(ad_h)^kg_1=g_2\},$$
    where we extend $g_1,g_2$ linearly to solutions in $\g\tensor\m_R\tensor\Omega_n$. Simplicial structure on $\HOM(g_1,g_2)$ is
    given by the one on $\{\Omega_n\}_{n\geq 0}$.
\end{definition}
From lemma \ref{Equivalences} it follows that $\DEL^\gamma(R)$ is indeed a simplicial groupoid.

Now consider in general a situation like in lemma \ref{Equivalences}: an $\li$-algebra $\g$, s.t. $\g=M\oplus\h$, where, if we forget the
differential, $\h$ is a Lie subalgebra. Suppose we have a morphism of $\li$-algebras $f:\g_1\rightarrow \g_2$, such that $f(M_1)\subset M_2$ and
$f(\h_1)\subset\h_2$.
\begin{lemma}\label{Partial}
    If $f$ is a quasi-isomorphism, it induces a weak equivalence between the simplicial Deligne groupoids, corresponding to $(\g_1,\h_1)$
    and $(\g_2,\h_2)$.
\end{lemma}
\textbf{Proof:}
    Applying cobar construction if necessary, we can assume that $\g_1,\g_2$ are dg Lie algebras. Suppose first that $f$ is an acyclic
    fibration. Then according to \cite{Hin1}(3.3.1), the
    corresponding functor $F$ between the Deligne groupoids of $\g_1$ and $\g_2$ is an acyclic fibration. That means in particular that
    $\pi_0(F)$ is an equivalence
    of categories and for any pair of objects $P,Q$ in $\DEL(\g_1)$, the map of simplicial sets $\HOM(P,Q)\rightarrow \HOM(F(P),F(Q))$ is
    an acyclic fibration. The simplicial set $\HOM(P,Q)$ has a subset $\HOM_{\h_1}(P,Q)$, consisting of the maps, defined by elements of
    $\h_1$. Since $f$ maps $\h_1$ to $\h_2$ and $M_1$ to $M_2$, it is clear that $F$ maps $\HOM_{\h_1}(P,Q)$ to $\HOM_{\h_2}(F(P),F(Q))$, and
    moreover the inverse image under $F$ of $\HOM_{\h_2}(F(P),F(Q))$ is $\HOM_{\h_1}(P,Q)$. Therefore, since pullback of an acyclic fibration is
    again an acyclic fibration, we conclude that $F$ defines an acyclic fibration from $\HOM_{\h_1}(P,Q)$ to $\HOM_{\h_2}(F(P),F(Q))$, and these are the
    mapping spaces in $\DEL(\g_1,\h_1)$ and $\DEL(\g_2,\h_2)$.

    If $f$ is not an acyclic fibration, it splits $f:g_1\overset{i}{\rightarrow} \g_3\overset{p}{\rightarrow} \g_2$, where $p$ is an acyclic fibration
    and the $i$ is an acyclic cofibration. Let $\h_3,M_3$ be the inverse images of $\h_2,M_2$ under $p$. Clearly $i^{-1}(\h_3)=\h_1$. From the
    construction of $\g_3$ (see e.g. \cite{Hin3}(2.2.4)) it is clear that one can define a quasi-isomorphism $\g_3\rightarrow \g_1$, left inverse
    to $i$, s.t. the image of $\h_3$ is in $\h_1$ and of $M_3$ is in $M_1$ (indeed just send all joined boundaries and coboundaries to $0$).
    Obviously this quasi-isomorphism is an acyclic fibration.
$\blacksquare$
\begin{proposition}\label{PartialA}
    Let $\gamma:A\rightarrow B$ and $\gamma':A'\rightarrow B'$ be morphisms of associative algebras.
    Let $q:A\rightarrow A',p:B\rightarrow B'$ be a quasi-isomorphism from $\gamma$ to $\gamma'$.
    Then the Deligne groupoids $\DEL^\gamma$ and $\DEL^{\gamma'}$ are weakly equivalent.
\end{proposition}
\textbf{Proof:}
    The map $p$ defines two maps: $C^\bullet(B',B')\rightarrow C^\bullet(B,B')$ and $C^\bullet(B,B)\rightarrow C^\bullet(B,B')$. Let $\g$ be the
    corresponding fiber product of $C^\bullet(B',B'),C^\bullet(B,B)$ over $C^\bullet(B,B')$ (as vector spaces). We can identify $\g$ with the
    subspace of $C^\bullet(B',B')\times C^\bullet(B,B)$, consisting of pairs of elements, whose images in $C^\bullet(B,B')$ coincide. Componentwise
    operations define of $\g$ the structure of a dg Lie algebra. We have canonical $f:\g\rightarrow C^\bullet(B,B)$ and $f':\g\rightarrow
    C^\bullet(B',B')$, and since
    $p$ is a quasi-isomorphism, these maps are quasi-isomorphisms. We extend by the means of $f,f'$ the $\li$ structures from
    $C^\bullet(A,A)\oplus s^{-1}C^\bullet(A,B)\oplus C^\bullet(B,B)$ and $C^\bullet(A,A)\oplus s^{-1}C^\bullet(A,B')\oplus C^\bullet(B',B')$ to
    $C^\bullet(A,A)\oplus s^{-1}C^\bullet(A,B)\oplus \g$ and $C^\bullet(A,A)\oplus s^{-1}C^\bullet(A,B')\oplus\g$ respectively. We have three
    quasi-isomorphisms:
        $$f:C^\bullet(A,A)\oplus s^{-1}C^\bullet(A,B)\oplus \g\rightarrow C^\bullet(A,A)\oplus s^{-1}C^\bullet(A,B)\oplus C^\bullet(B,B),$$
        $$f':C^\bullet(A,A)\oplus s^{-1}C^\bullet(A,B')\oplus \g\rightarrow C^\bullet(A,A)\oplus s^{-1}C^\bullet(A,B')\oplus C^\bullet(B',B'),$$
        $$p_*:C^\bullet(A,A)\oplus s^{-1}C^\bullet(A,B)\oplus\g\rightarrow C^\bullet(A,A)\oplus s^{-1}C^\bullet(A,B')\oplus\g.$$
    So by lemma \ref{Partial}, the simplicial Deligne groupoids that correspond to $\gamma$ and $p\gamma$ are weakly equivalent. Doing the same
    thing with $q$ we get the final result.
$\blacksquare$

\section{Solution of the deformation problem in case of non-positively graded algebras}\label{Solution}

In this section we prove that if $\gamma:A\rightarrow B$ is a morphism of non-positively graded algebras, then $\DEF^\gamma$ is weakly
equivalent to $\DEL^\gamma$. We make this requirement because all non-positively graded almost free algebras are cofibrant, whereas in general
it is not true for $\mathbb Z$-graded algebras.

Let $(M,\delta)$ be a differential $\Z$-graded associative $R$-algebra. The bar construction of $(M,\delta)$ is the codifferential coassociative
coalgebra $(\B M,\B\delta)$, where $\B M:=T(sM)$ and $\B\delta=s\circ(\delta+\mu)\circ s^{-1}$, where $\mu$ is the multiplication on $M$. In
turn for a codifferential coassociative coalgebra $(Z,\delta)$ the co-bar construction is the dg associative algebra $(\Omega Z,\Omega\delta)$,
where $\Omega Z:=T(s^{-1}Z)$, $\Omega\delta:=s^{-1}\circ(\delta+\Delta)\circ s$, where $\Delta$ is the comultiplication on $Z$. We will denote
the co-bar construction on the bar construction of $M$ by $\Omega\B M$.

$\Omega\B M$ is an almost free non-positively graded algebra, and it is cofibrant if $M$ is non-positively graded. Indeed almost free algebras
are cofibrant in the category of non-positively graded dg associative algebras (\cite{Get}(4.6)), and an acyclic fibration of $\mathbb Z$-graded
algebras induces an acyclic fibration of their truncations at 0, with 0-part consisting of cocycles. Left lifting property then goes over to the
category of all $\mathbb Z$-graded algebras. There is a natural transformation $\epsilon:\Omega\B\rightarrow Id_{\DAlg(R)}$.

Clearly $\Omega\B$ is a functor, hence it extends to $\Mor\DAlg(R)$, and we will denote the extension again by $\Omega\B$. If $M$ and $N$ in
$\phi:M\rightarrow N\in\Mor\DAlg(R)$ are non-positively graded algebras, then $\Omega\B M,\Omega\B N$ are cofibrant, and if in addition $\phi$
is a cofibration, then $\Omega\B\phi$ is a cofibration. Indeed, if $\phi$ is injective then, since $\Omega\B\phi$ maps generators of the domain
injectively to generators of the co-domain, it is a cofibration (\cite{Get} page 42), if $\phi$ is not injective we can split it into an
injective cofibration, followed by an acyclic fibration, and then $\Omega\B\phi$ is a retract of a cofibration. Natural transformation
$\epsilon$ extends to $\Omega\B\rightarrow Id_{\Mor\DAlg(R)}$.

Let $C(R)$ be the simplicial subcategory of $\Mor\DALG(R)$, consisting of cofibrations between non-positively graded, cofibrant algebras. As
noted above $\Omega\B$ maps $C(R)$ into itself, and there is a natural transformation $\epsilon:\Omega\B\rightarrow Id_{C(R)}$.
\begin{lemma}\label{Bar}
    $\Omega\B:C(R)\rightarrow \Omega\B(C(R))$ and the inclusion of $\Omega\B(C(R))$ in $C(R)$ induce a weak equivalence
    of the nerves of these two categories.
\end{lemma}
Now consider a morphism $\gamma:A\rightarrow B\in\Mor\DAlg(\K)$ between non-positively graded algebras that we want to deform. By proposition
\ref{PartialA}, for the purpose of comparing $\DEL^\gamma$ with $\DEF^\gamma$, we can consider $\gamma$ as a cofibration between cofibrant
objects, i.e. $\gamma\in C(R)$.

We have the $\li$-algebra $(\g,d^g)$ that describes deformations of $\gamma$ as an $\A$-structure. These deformations consist of codifferentials
on $\B(A\tensor R),\B(B\tensor R)$ and coalgebra morphisms $\B(A\tensor R)\rightarrow\B(B\tensor R)$. Applying $\Omega$ to them we get objects
in $C(R)$, whose reduction modulo $\m_R$ is $\Omega\B\gamma$. In this way for every solution $g'$ of MCE in $(\g\tensor\m_R,d^g)$ there is a
corresponding object $\F(g')$ in $\DEF^\gamma$. Moreover, mapping spaces in the Deligne groupoid represent morphisms between morphisms of
coalgebras (lemma \ref{Equivalences}), therefore since the cobar construction is a functor, $\F$ is actually a simplicial functor
$\DEL^\gamma(R)\rightarrow\DEF^\gamma(R)$. Clearly $\Omega$ is functorial in $R$, hence $\underline{F}$ is a morphism in $Hom(\dgart(\K),\GRP)$.
\begin{theorem}
    $\F$ is a weak equivalence $\DEL^\gamma(R)\rightarrow\DEF^\gamma(R)$.
\end{theorem}
\textbf{Proof:}
    As noted above, $\F$ maps $\DEL^\gamma(R)$ to $C(R)$, and the latter is a full simplicial subcategory of the category of cofibrant objects
    in $\Mor\DALG(R)$. Clearly the image of $\F$ is in the fiber of $\PI_R$ over $\BAMMA$ (see definition \ref{SecondDefinition}).
    To prove that $\F$ is a weak equivalence of $\DEL^\gamma(R)$ with this fiber note that $\F$ maps $\DEL^\gamma(R)$ identically on the
    fiber of $\Omega\B(C(R))\rightarrow\Omega\B(C(\K))$ at $\BAMMA$. To compare this fiber with the one of $\PI_R$, note that
    a map between simplicial groupoids is a weak equivalence if and only if it induces a weak equivalence of the nerves. Moreover, since for
    fibrations between simplicial categories, the nerve of a fiber is equivalent to the fiber of the nerve (see proof of lemma
    \ref{Fiber}), from lemma \ref{Bar} we conclude that the fiber of $\Omega\B(C(R))\rightarrow\Omega\B(C(\K))$ at $\BAMMA$ is
    indeed weakly equivalent to the fiber of $\PI_R$ at $\BAMMA$.
$\blacksquare$

\subsection{Cohomology}

If we start with a pair of associative algebras $A,B$, concentrated in degree 0, and an associative algebra morphism $\gamma:A\rightarrow B$,
then the differential $d^g_1$ on $\g$ is as follows ($\alpha\in C^\bullet(A,A)$)
    $$d^g_1(\alpha)(a_1\tensor...\tensor a_{n+1})=(-1)^{n+1}(-a_1\alpha(a_1\tensor...\tensor a_{n+1})+$$
    $$+\underset{1\leq i\leq n}{\Sigma}(-1)^{i+1}\alpha(a_1\tensor...\tensor a_ia_{i+1}\tensor...\tensor
    a_{n+1})+(-1)^n\alpha(a_1\tensor...\tensor a_n)a_{n+1})-$$
    $$-\gamma(\alpha(a_1\tensor...\tensor a_n)).$$
In other words
    $$d^g_1(\alpha)=(-1)^{\alpha+1}HD(\alpha)-\gamma_*(\alpha),$$
where $HD$ is the Hochschild differential and the degree of $\alpha$ is the number of its arguments minus 1. For a $\gamma'\in C^\bullet(A,B)$
we have
    $$d^g_1(\gamma')(a_1\tensor...\tensor a_{n+1})=\gamma'(a_1\tensor...\tensor a_n)\gamma(a_{n+1})+$$
    $$+(-1)^{n+1}\gamma(a_1)\gamma'(a_2\tensor...\tensor a_{n+1})+(-1)^n\underset{1\leq i\leq
    n}{\Sigma}(-1)^{i+1}\gamma'(a_1\tensor...\tensor a_ia_{i+1}\tensor...\tensor a_{n+1}),$$
that is
    $$d^g_1(\gamma')=(-1)^{\gamma'} HD(\gamma').$$
The case of  $\beta\in C^\bullet(B,B)$ is similar to $\alpha\in C^\bullet(A,A)$. In total we can describe the differential complex
    $$(C^\bullet(A,A)\oplus C^\bullet(B,B)\oplus C^\bullet(A,B),d^g_1)$$
as the cone of the morphism of complexes:
    $$(C^\bullet(A,A),(-1)^nHD)\oplus(C^\bullet(B,B),(-1)^nHD)\rightarrow(C^\bullet(A,B),(-1)^nHD),$$
where $\alpha+\beta$ is mapped to $\gamma^*(\beta)-\gamma_*(\alpha)$. Gerstenhaber and Schack (\cite{GS2} page 249, \cite{GS3} page 8) have
constructed cohomology of a morphism using the same complexes with untwisted Hochschild differentials and the opposite morphism of complexes.
The cohomology in their construction is the same as in ours.

Cohomology of an $\li$-algebra has a canonical structure of a Lie algebra (induced by the binary operation of the $\li$-structure). The Lie
structure on cohomology appearing in \cite{GS2},\cite{GS3} is induced by the $\li$-structure on cochains that we have described, indeed the
binary operation of this $\li$-structure coincides with the commutator of the circle operation on cochains that is used in
\cite{GS2},\cite{GS3}.

\end{document}